\documentclass[12pt,a4paper]{amsart}
\usepackage{amsmath, amsthm, amsfonts, amssymb, graphicx, hyperref, bm,verbatim,siunitx}

\theoremstyle{plain}
\newtheorem{thrm}{Theorem}
\newtheorem{lmm}[thrm]{Lemma}

\newtheorem{cnjctr}{Conjecture}
\newtheorem{prblm}{Problem}

\numberwithin{sblmm}{thrm} 
\numberwithin{equation}{section}

\renewcommand{\phi}{\varphi}

\newcommand{\Mod}[1]{\ (\mathrm{mod}\ #1)}

\begin{document}
\title{Gaps between primes}
\author{James Maynard}
\address{Mathematical Institute, Oxford, England OX1 4AU}
\email{james.alexander.maynard@gmail.com}
\begin{abstract}
We discuss recent advances on weak forms of the Prime $k$-tuple Conjecture, and its role in proving new estimates for the existence of small gaps between primes and the existence of large gaps between primes.
\end{abstract}
\maketitle

\section{Introduction}
It follows from the Prime Number Theorem that the \textit{average} gap between primes less than $X$ is of size roughly $\log{X}$ when $X$ is large. We expect, however, that occasionally these gaps are rather smaller than $\log{X}$, and occasionally rather larger. Specifically, based on random models and numerical evidence, we believe that the largest and smallest gaps are as described in the following two famous conjectures\footnote{The original conjecture of Cram\'er \cite{Cramer} was a stronger statement which is no longer fully believed to be true. We expect the weaker version given here to hold.}.
\begin{cnjctr}[Twin Prime Conjecture]\label{cnjctr:TwinPrime}
There are infinitely many pairs of primes which differ by exactly 2.
\end{cnjctr}

\begin{cnjctr}[Cram\'er's Conjecture, weak form]\label{cnjctr:Cramer} Let $p_n$ denote the $n^{th}$ prime. Then
\[
\sup_{p_n\le X}(p_{n+1}-p_n)=(\log{X})^{2+o(1)}.
\]
\end{cnjctr}
Moreover, the Twin Prime Conjecture can be thought of as a special case of the far-reaching Prime $k$-tuple Conjecture describing more general patterns of many primes.
\begin{cnjctr}[Prime $k$-tuple Conjecture]\label{cnjctr:KTuples}
Let $L_1,\dots,L_k$ be integral linear functions $L_i(n)=a_in+b_i$ such that for every prime $p$ there is an integer $n_p$ with $\prod_{i=1}^k L_i(n_p)$ coprime to $p$. Then there are infinitely many integers $n$ such that all of $L_1(n),\dots,L_k(n)$ are primes.
\end{cnjctr}
%We say that $\{L_1,\dots,L_k\}$ are \textit{admissible} if 

Unfortunately all these conjectures seem well beyond the current techniques, but we are able to make partial progress by showing that we do have gaps which are smaller or larger than the average gap. A key part of recent progress on results about gaps between primes has been weak versions of Conjecture \ref{cnjctr:KTuples}, where one shows there are infinitely many integers $n$ such that \textit{several} (rather \textit{all}) of the linear functions are prime at $n$.
\begin{thrm}[\cite{Maynard}]\label{thrm:WeakTuples}
Let $L_1,\dots,L_k$ be integral linear functions $L_i(n)=a_in+b_i$ such that for every prime $p$ there is an integer $n_p$ with $\prod_{i=1}^k L_i(n_p)$ coprime to $p$. Then there is a constant $c>0$ such that there are infinitely many integers $n$ where at least $c\log{k}$ of $L_1(n),\dots,L_k(n)$ are primes. 
\end{thrm}
Variants of Theorem \ref{thrm:WeakTuples} have been important in recent results on small and large gaps between primes, and have proven useful because the method of proof is quite flexible and can generalize to other situations. Specifically, we now know the following approximations to Conjecture \ref{cnjctr:TwinPrime} and Conjecture \ref{cnjctr:Cramer}.
\begin{thrm}[\cite{Polymath}]\label{thrm:SmallGaps}
There are infinitely many pairs of primes which differ by at most $246$.
\end{thrm}
 \begin{thrm}[\cite{FGKMT}]\label{thrm:LargeGaps}
 There exists a constant $c>0$ such that
 \[
 \sup_{p_n\le X}(p_{n+1}-p_n)\ge c \frac{\log{X}\cdot\log\log{X}\cdot \log\log\log\log{X}}{\log\log\log{X}} .
 \]
 \end{thrm}

Moreover, in the case of small gaps between primes, we can show the existence of \textit{many} primes in bounded length intervals.
\begin{thrm}[\cite{MaynardDense}]\label{thrm:ManyGaps}
There exists a constant $c>0$ such that, for all $X\ge 2$ there are at least $c X\exp(-\sqrt{\log{X}})$ integers $x\in[X,2X]$ such that
\[
\#\{\text{primes in }[x,x+y]\} \ge c\log{y}.
\]
\end{thrm}
For example, this shows that there are infinitely many intervals of length $e^{m/c}$ containing $m$ primes, and, for fixed $\epsilon>0$, infinitely many $x$ such that the interval $[x,x+(\log{x})^\epsilon]$ contains $\epsilon c\log\log{x}$ primes.

\section{The GPY sieve method}
We aim to prove Theorem \ref{thrm:WeakTuples} by the `GPY method', which can be interpreted as a first moment method which was introduced to study small gaps between primes by Goldston, Pintz and Y\i ld\i r\i m \cite{GPY}. This takes the following basic steps, for some given set $\{L_1,\dots,L_k\}$ of distinct functions satisfying the hypotheses of Conjecture \ref{cnjctr:KTuples}:
\begin{enumerate}
\item We choose a probability measure $w$ supported on integers in $[X,2X]$.
\item We calculate the expected number of the functions $L_i(n)$ which are prime at $n$, if $n$ is chosen randomly with probability $w(n)$.
\item If this expectation is at least $m$, then there must be some $n\in[X,2X]$ such that at least $m$ of $L_1(n),\dots, L_k(n)$ are primes.
\item If this holds for all large $X$, then there are infinitely many such integers $n$.
\end{enumerate}
This procedure only works if we can find a probability measure $w$ which is suitably concentrated on integers $n$ where many of the $L_i(n)$ are prime, but at the same time is simple enough that we can calculate this expectation unconditionally. We note that by linearity of expectation, it suffices to be able to calculate the probability that any one of the linear functions is prime. However, any slowly changing smooth function $w$ is insufficient since the primes have density 0 in the integers, whereas any choice of $w$ explicitly depending on the joint distribution of prime values of the $L_i$ is likely to be too complicated to handle unconditionally. 

Sieve methods are a flexible set of tools, developed over the past century, which provide natural choices for the probability measure $w$. The situation of simultaneous prime values of $L_1,\dots,L_k$ is a `$k$-dimensional' sieve problem. For such problems when $k$ is large but fixed, the Selberg sieve tends to be the type of sieve which performs best. The standard choice of Selberg sieve weights (which are essentially optimal in closely related situations) are
\begin{equation}
w(n)\propto\Bigl(\sum_{\substack{d|\prod_{i=1}^k L_i(n) \\ d<R}}\mu(d)\Bigl(\log\frac{R}{d}\Bigr)^k\Bigr)^2,\label{eq:SelbergWeights}
\end{equation}
where $R$ is a parameter which controls the complexity of the sieve weights, and $w$ is normalized to sum to $1$ on $[X,2X]$.

To calculate the probability that $L_j(n)$ is prime with this choice of $w(n)$, we wish to estimate the sum of $w(n)$ over $n\in[X,2X]$ such that $L_j(n)$ is prime. To do this we typically expand the divisor sum in the definition \eqref{eq:SelbergWeights} and swap the order of summation. This reduces the problem to estimating the number of prime values of $L_j(n)$ for $n\in [X,2X]$ in many different arithmetic progressions with moduli of size about $R^2$. The Elliott-Halberstam conjecture \cite{ElliottHalberstam} asserts that we should be able to do this when $R^2<X^{1-\epsilon}$, but unconditionally we only know how to do this when $R^2<X^{1/2-\epsilon}$, using the Bombieri-Vinogradov Theorem \cite{Bombieri,Vinogradov}. After some computation one finds that, provided we do have suitable estimates for primes in arithmetic progressions, the choice \eqref{eq:SelbergWeights} gives
\begin{equation}
\mathbb{E}\,\#\{1\le i\le k:L_i(n)\text{ prime}\}=\Bigl(2-\frac{2}{k+1}+o(1)\Bigr)\frac{\log{R}}{\log{X}}. \label{eq:SelbergExpectation}
\end{equation}
In particular, this is less than $1$ for all large $X$, even if we optimistically assume the Elliott-Halberstam conjecture and take $R\approx X^{1/2-\epsilon}$, and so it appears that we will not be able to conclude anything about primes in this manner.

The groundbreaking work of Goldston-Pintz-Y\i ld\i r\i m \cite{GPY} showed that a variant of this choice of weight actually performs much better. They considered
\begin{equation}
w(n)\propto\Bigl(\sum_{\substack{d|\prod_{i=1}^k L_i(n) \\ d<R}}\mu(d)\Bigl(\log\frac{R}{d}\Bigr)^{k+\ell}\Bigr)^2,\label{eq:GPYWeights}
\end{equation}
where $\ell$ is an additional parameter to be optimized over. With the choice \eqref{eq:GPYWeights}, one finds that provided we have suitable estimates for primes in arithmetic progressions, we obtain
\[
\mathbb{E}\,\#\{1\le i\le k:L_i(n)\text{ prime}\}=\Bigl(4+O\Bigl(\frac{1}{\ell}\Bigr)+O\Bigl(\frac{\ell}{k}\Bigr)\Bigr)\frac{\log{R}}{\log{X}}.
\]
This improves upon \eqref{eq:SelbergExpectation} by a factor of about 2 when $k$ is large and $\ell\approx k^{1/2}$. This falls just short of proving that two of the linear functions are simultaneously prime when $R=X^{1/4-\epsilon}$ as allowed by the Bombieri-Vinogradov theorem, but any small improvement allowing $R=X^{1/4+\epsilon}$ would give bounded gaps between primes! By considering additional possible primes, Goldston Pintz and Y\i ld\i r\i m were able to show
\[
\liminf_n \frac{p_{n+1}-p_n}{\log{p_n}}=0,
\]
 finally extending a long sequence of improvements to upper bounds for the left hand side \cite{Erdos, Rankin, Ricci,BombieriDavenport,Piltai,Uchiyama,Huxley,Huxley2,Huxley3,Maier}.
 
   Building on earlier work of  Bombieri, Fouvry, Friedlander and Iwaniec \cite{FouvryIwaniec,BFIi,BFIii,BFIiii}, Zhang \cite{Zhang} succeeded in establishing an extended version of the Bombieri-Vinogradov Theorem for moduli with no large prime factors allowing $R=X^{1/4+\epsilon}$, and ultimately this allowed him to show
\[
\mathbb{E}\,\#\{1\le i\le k:L_i(n)\text{ prime}\}>1,
\]
if $k>\num{3500000}$ and $X$ is sufficiently large. The key breakthrough in Zhang's work was this result on primes in arithmetic progressions of modulus slightly larger than $X^{1/2}$. By choosing suitable linear functions, this then showed that
\[
\liminf_n(p_{n+1}-p_n)\le 7\cdot 10^7.
\]

\section{A modified GPY sieve method}\label{sec:Sieve}
An alternative approach to extending the work of Goldston, Pintz and Y\i ld\i r\i m  was developed independently by the author \cite{Maynard} and Tao (unpublished). The key difference was to consider the multidimensional generalization%\footnote{Selberg had Goldston and Y\i ld\i r\i m also had earlier work which considered similar setups, which were not quite sufficient for bounded gaps between primes}
\begin{equation}
w(n)\propto\Bigl(\sum_{\substack{d_1,\dots,d_k \\ d_i|L_i(n) \\ \prod_{i=1}^k d_i<R}}\mu(d)F\Bigl(\frac{\log{d_1}}{\log{R}},\dots,\frac{\log{d_k}}{\log{R}}\Bigr)\Bigr)^2,\label{eq:NewWeights}
\end{equation}
for suitable smooth functions $F:\mathbb{R}^k\rightarrow \mathbb{R}$ supported on $[0,\infty)^k$. The flexibility of allowing the function $F$ to depend on each divisor $d_1,\dots,d_k$ of $L_1(n),\dots,L_k(n)$ allows us to make $w(n)$ more concentrated on integers $n$ when many of the $L_i(n)$ are prime.

With this choice, after some computation, one finds that 
\begin{equation}
\mathbb{E}\,\#\{1\le i\le k:L_i(n)\text{ prime}\}=\Bigl(\frac{\sum_{i=1}^k J_i(\tilde{F})}{I(\tilde{F})}+o(1)\Bigr)\frac{\log{R}}{\log{X}},
\label{eq:NewExpectation}
\end{equation}
provided, as before, we are able to count primes in arithmetic progressions to modulus $R^2$ on average. Here the $J_i(\tilde{F})$ and $I(\tilde{F})$ are $k$-dimensional integrals depending on a transform\footnote{$\tilde{F}$ is $F$ differentiated with respect to each coordinate.} $\tilde{F}$ of $F$, given by
\begin{align*}
J_\ell(\tilde{F}) &=\idotsint\limits_{\sum_{i\ne \ell} t_i\le 1}\Bigl(\int_0^{1-\sum_{i\ne \ell}t_i} \tilde{F}(t_1,\dots,t_k)dt_\ell \Bigr)^2 dt_1\dots dt_{\ell-1}dt_{\ell+1}\dots dt_k,\\
I(\tilde{F}) &=\idotsint\limits_{\sum_{i=1}^k t_i\le 1}\tilde{F}(t_1,\dots,t_k)^2dt_1\dots dt_k.
\end{align*}
For any piecewise smooth choice of $\tilde{F}$ supported on $\sum_{i=1}^kt_i\le 1$ there is a corresponding choice of $F$. In particular, we can show that many of the $L_i$ are simultaneously prime infinitely often, if we can show that $\sup_{\tilde{F}}(\sum_{i=1}^k J_i(\tilde{F})/I(\tilde{F}))\rightarrow \infty$ as $k\rightarrow\infty$.

A key advantage of this generalization is that we can make use of high dimensional phenomena such as concentration of measure. We concentrate on functions $\tilde{F}$ the form
\[
\tilde{F}(t_1,\dots,t_k)=\begin{cases}
\prod_{i=1}^k G(t_i),\qquad &\sum_{i=1}^kt_i<1,\\
0,&\text{otherwise,}
\end{cases}
\]
and make a probabilistic interpretation of the integrals $J_i(\tilde{F})$ and $I(\tilde{F})$. Let $Z_1,\dots,Z_k$ be i.i.d. random variables on $[0,\infty]$ with probability density function $G^2$, expectation $\mu=\int_0^\infty t G(t)^2dt$ and variance $\sigma^2=\int_0^\infty (t-\mu)^2G(t)^2dt$. Then
\begin{align*}
I(\tilde{F})&=\mathbb{P}\Bigl(\sum_{i=1}^k Z_i<1\Bigr),\qquad
J_\ell(\tilde{F})\ge \Bigl(\int_0^{1/2}G(t)dt\Bigr)^2\mathbb{P}\Bigl(\sum_{i=1}^k Z_i<\frac{1}{2}\Bigr).
\end{align*}
The random variable $\sum_{i=1}^k Z_i$ has mean $k\mu$ and variance $k\sigma^2$, and so it becomes concentrated on $k\mu$ when $k$ is large provided $\sigma^2/\mu\rightarrow 0$ as $k\rightarrow\infty$. In particular, if $\mu<1/3k$ and $\sigma^2k\rightarrow 0$, then $\mathbb{P}(\sum_{i=1}^k Z_i<1/2)$ approaches 1 as $k\rightarrow\infty$. Therefore, to show that $\sum_{\ell=1}^kJ_\ell(\tilde{F})/I(\tilde{F})\rightarrow\infty$ as $k\rightarrow\infty$, it suffices to find a function $G$ satisfying
\[
\int_0^\infty t G(t)^2<\frac{1}{3k},\qquad\int_0^\infty G(t)^2 dt=1,
\]
\[
k\int_0^\infty t^2G(t)^2dt\rightarrow 0 \quad \text{and}\quad k\Bigl(\int_0^\infty G(t)dt\Bigr)^2\rightarrow \infty \quad \text{as}\quad k\rightarrow \infty.
\]
We find that choosing
\begin{equation}
G(t)\approx \begin{cases}
\frac{\sqrt{k\log{k}}}{1+t k\log{k}},\qquad &t<k^{-3/4},\\
0,&\text{otherwise,}
\end{cases}\label{eq:SmoothChoice}
\end{equation}
gives a function $G$ satisfying these constraints. (This choice can be found via the calculus of variations, and further calculations show that this choice is essentially optimal.) Putting this all together, we find that for some constant $c>0$
\[
\mathbb{E}\,\#\{i:\,L_i(n)\text{ prime}\}\ge ( c\log{k}+o(1))\frac{\log{R}}{\log{X}}.
\]
In particular, taking $R=X^{1/4-o(1)}$ (as allowed by the Bombieri-Vinogradov theorem) and letting $k$ be sufficiently large, we find that there are infinitely many integers $n$ such that $(c\log{k})/4+o(1)$ of the $L_i(n)$ are simultaneously prime. Performing these calculations carefully allows one to take $c\approx 1$ when $k$ is large.

Morally, the effect of such a choice of function $\tilde{F}$ is to make `typical' divisors $(d_1,\dots,d_k)$ occurring in \eqref{eq:NewWeights} to have $\prod_{i=1}^k d_i$ smaller, but for it to be more common that some of the components of $(d_1,\dots,d_k)$ are unusually large when compared with \eqref{eq:SelbergWeights} or \eqref{eq:GPYWeights}. This correspondingly causes the random integer $n$ chosen with probability $w(n)$ to be such that the $L_i(n)$ are more likely to have slightly smaller prime factors, but it is also more likely that some of the $L_i(n)$ are prime.

\section{Consequences for small gaps between primes}

The above argument shows that there is a constant $c>0$ such that for any set $\{L_1,\dots,L_k\}$ of integral linear functions satisfying the hypotheses of Conjecture \ref{cnjctr:KTuples}, at least $c\log{k}$ of the linear functions $L_i$ are simultaneously prime infinitely often. To show that there are primes close together, we simply take the linear functions to be of the form $L_i(n)=n+h_i$ for some integers $h_1\le \dots \le h_k$ chosen to make $\prod_{i=1}^k L_i(n)$ have no fixed prime divisor, and so that $h_k-h_1$ is small. 

In general, a good choice of the $h_i$ is to take $h_i$ to be the $i^{th}$ prime after the integer $k$. With this choice, $\prod_{i=1}^kL_i(n)$ has no fixed prime divisor and $h_k-h_i\approx k\log{k}$, so we can find intervals of length $k\log{k}$ containing $c\log{k}$ primes infinitely often. By working out the best possible implied constants, the argument we have sketched allows us to show that there are $m$ primes in an interval of length $O(m^3e^{4m})$ infinitely often. By incorporating refinements of the work of Zhang on primes in arithmetic progressions by the Polymath 8a project \cite{Polymath8A}, and using some bounds from Harman's sieve \cite{HarmanBook}, this can be improved slightly to $O(e^{3.815m})$ by work of Baker-Irving \cite{BakerIrving}.

If we are only interested in just how small a \textit{single} gap can be, then we can improve the analysis and get an explicit bound on the size of the gap by adopting a numerical analysis perspective. The key issue is to find the smallest value of $k$ such that for all large $X$
\[
\mathbb{E}\,\#\{1\le i\le k:\,L_i(n)\text{ prime}\}>1,
\]
since this immediately implies that two of the linear functions are simultaneously prime. Recalling that we can choose $R=x^{1/4-\epsilon}$ for any $\epsilon>0$, using \eqref{eq:NewExpectation} we reduce the problem to finding a value of $k$ as small as possible, such that we can find a function $\tilde{F}:[0,\infty)^k\rightarrow \mathbb{R}$ with
\[
\frac{\sum_{i=1}^k J_i(\tilde{F})}{I(\tilde{F})}>4.
\]
We fix some basis functions $g_1,\dots,g_r:[0,\infty)^k\rightarrow \mathbb{R}$ which are supported on $(t_1,\dots,t_k)$ satisfying $\sum_{i=1}^k t_i\le 1$, and restrict our attention to functions $\tilde{F}$ in the linear span of the $g_i$; that is functions $\tilde{F}$ of the form
\[
\tilde{F}(t_1,\dots,t_k)=\sum_{i=1}^r f_i g_i(t_1,\dots,t_r),
\]
for some coefficients $\mathbf{f}=(f_1,\dots f_r)\in\mathbb{R}^r$ which we think of as variables we will optimize over. For such a choice, we find that $\sum_{i=1}^k J_i(\tilde{F})$ and $I(\tilde{F})$ are both quadratic forms in the variables $f_1,\dots,f_r$, and the coefficients of these quadratic forms are explicit integrals in terms of $g_1,\dots,g_r$. If we choose a suitably nice basis $g_1,\dots,g_r$, then these integrals are explicitly computable, and so we obtain explicit $r\times r$ real symmetric matrices $M_1$ and $M_2$ such that
\[
\frac{\sum_{i=1}^k J_i(\tilde{F})}{I(\tilde{F})}=\frac{\mathbf{f}^T M_1\mathbf{f}}{\mathbf{f}^T M_2\mathbf{f}}.
\]
We then find that the choice of coefficients $\mathbf{f}$ which maximizes this ratio is the eigenvector of $M_2^{-1}M_1$ corresponding to the largest eigenvalue, and the value of the ratio is given by this largest eigenvalue. Thus the existence of a good function $\tilde{F}$ can be reduced to checking whether the largest eigenvalue of a finite matrix is larger than 4, which can be performed numerically by a computer.

  If we choose the $g_i$ to be symmetric polynomials of low degree, then this provides a nice basis since the corresponding integrals have a closed form solution, and allows one to make arbitrarily accurate numerical approximations to the optimal function $\tilde{F}$ with enough computation. For the problem at hand, these numerical calculations are large but feasible. This approach ultimately allows one to show that if $k=54$ there is a function $\tilde{F}$ such that $\sum_{i=1}^{54}J_i(\tilde{F})/I(\tilde{F})>4$. 
  
  To turn this into small gaps between primes we need to choose the shifts $h_i$ in our linear functions $L_i(n)=n+h_i$ so that $\prod_{i=1}^{54}L_i(n)$ has no fixed prime divisor, and the $h_i$ are in as short an interval as possible. This is a feasible exhaustive numerical optimization problem, with an optimal choice of the $\{h_1,\dots,h_{54}\}$ given by \footnote{Such computations were first performed by Engelsma - see \url{http://www.opertech.com/primes/k-tuples.html}}
  \begin{align*}
  \{&0, 2, 6, 12, 20, 26, 30, 32, 42, 56, 60, 62, 72, 74, 84, 86, 90, 96, \\
  &\quad  104, 110, 114, 116, 120, 126, 132, 134, 140, 144, 152, 156, \\
  &\quad 162, 170, 174, 176, 182, 186, 194, 200, 204, 210, 216, 222, \\
  &\quad  224, 230, 236, 240, 242, 246, 252, 254, 260, 264, 266, 270\}.
  \end{align*}
  Putting this together, we find
\[
\liminf_n (p_{n+1}-p_n)\le 270.
\]
This result is not quite the current record - in \cite{Polymath} we make some further technical refinements to the sieve, which corresponds to modifying the expressions $J_i(\tilde{F})$ and $I(\tilde{F})$ slightly. This ultimately allows us to improve $k=54$ to $k=50$ and correspondingly improve from gaps of size at most 270 to gaps of at most 246, giving Theorem \ref{thrm:SmallGaps}. The main ideas are the same as above.

\section{Large gaps between primes}

It turns out that because Theorem \ref{thrm:WeakTuples} gives strong partial information about the joint distribution of prime solutions to linear equations, which lie at the heart of many basic estimates about primes, it can be used to make progress on the existence of \textit{large} gaps between primes, although the connection is less direct than with small gaps.

The first major breakthrough on large gaps between primes was due to Westzynthius \cite{Westzynthius}, who showed that there were gaps which could be arbitrarily large compared with the average gap. During the 1930s this was refined with ideas due to Erd\H os \cite{Erdos} and Rankin \cite{RankinLarge}, giving
\begin{equation}
 \sup_{p_n\le X}(p_{n+1}-p_n)\ge c' \frac{\log{X}\cdot\log\log{X}\cdot \log\log\log\log{X}}{(\log\log\log{X})^2},\label{eq:RankinBound}
\end{equation}
for some constant $c'>0$. Subsequent improvements over the next 75 years \cite{Schonhage,Rankin2,MaierPomerance} were only in improving the value of the constant $c'$, the strongest being $c'=2e^\gamma+o(1)$, due to Pintz \cite{Pintz}. All of these approaches are based on a variant of the following lemma, which reduces the problem of constructing large gaps between primes to a combinatorial covering problem.
\begin{lmm}\label{lmm:CRT}
If one can choose residue classes $a_p\Mod{p}$ for $p\le x$ such that every element of $\{1,\dots,y\}$ is congruent to $a_p\Mod{p}$ for some $p\le x$, then there is a prime $p_n\ll e^x$ such that $p_{n+1}-p_n\ge y$.
\end{lmm}
This lemma (which is a simple consequence of the Chinese Remainder Theorem), can be thought of as a natural generalization of the argument that $n!+2,\dots,n!+n$ are $n$ consecutive composite integers, and so explicitly demonstrates a prime gap of size at least $n$. (This roughly corresponds to choosing $a_p=1$ for all primes $p$ and taking $y=x$.)

The key idea in the Erd\H os-Rankin construction was to choose $a_p=0$ for `medium sized' primes, and choose $a_p$ differently for small and large primes. Specifically, a version of their argument follows the following strategy to choose the $a_p$ in turn, for some parameter $2<z<x^{1/2}$: 
\begin{enumerate}
\item Choose $a_p=0$ for `medium primes' $p\in [z,x/3]$.
\item Choose $a_p=1$ for `small primes' $p<z$.
\item Choose $a_p$ greedily for `large primes' $p\in (x/3,x]$.
\end{enumerate}
By `choosing greedily' we mean that we pick a residue class $a_q\Mod{q}$ which contains the largest number of integers in $\{1,\dots,y\}$ which are not congruent to $a_p\Mod{p}$ for some previously chosen $a_p$. There must be a residue class containing at least one uncovered element if we have not already covered all elements.

By choosing $z$ appropriately\footnote{$z=\exp(\log{x}\cdot \log\log\log{x}\,/\,2\log\log{x})$ gives a suitable choice.}, this allows one to cover $\{1,\dots,y\}$ by residue classes $a_p\Mod{p}$ for $p\le 7y(\log\log{y})^2/(\log{y}\cdot\log\log\log{y})$, which ultimately results in the bound \eqref{eq:RankinBound}. The key feature making this construction work is that choosing $a_p=0$ for a very large number of `medium primes' is much more efficient that a typical choice. This is because integers in $\{1,\dots,y\}$ with no prime factors bigger than $z_1$ are much less common than integers avoiding a random residue class for each prime bigger than $z_1$.

It was a well-known challenge of Erd\H os as to whether one could improve upon \eqref{eq:RankinBound} by an arbitrarily large constant, and this was verified independently by Ford-Green-Konyagin-Tao \cite{FGKT} and the author \cite{MaynardLarge}. Ultimately the approach of Ford-Green-Konyagin-Tao relied on the work of Green-Tao and Green-Tao-Ziegler on linear equations in primes \cite{GreenTao,GreenTao2,GreenTao3}, whereas the work of the author relied on versions of Theorem \ref{thrm:WeakTuples}. So far this second approach has proved more flexible for gaining quantitative improvements over \eqref{eq:RankinBound}, the strongest known results being due to a collaboration between all these authors \cite{FGKMT}.% following the overall approach based on Theorem \ref{thrm:WeakTuples}, but incorporating simplifications from the work of \cite{FGKT} as well as newer combinatorial ideas.

We focus here on the approach based around Theorem \ref{thrm:WeakTuples}, first thinking about obtaining an arbitrarily large constant improvement over \eqref{eq:RankinBound}. We follow the same overall strategy as Erd\H os-Rankin, but improve the analysis for the large primes. By choosing the residue classes in a more sophisticated manner, we are able to remove \textit{many} uncovered elements on average, rather than just 1 element, and this is the key feature which allows us to improve on \eqref{eq:RankinBound}.

Using the same choice of $a_p$ as Erd\H os-Rankin for $p\le x/3$, we see that the elements of $\{1,\dots,y\}$ which are not covered by $a_p\Mod{p}$ for $p\le x/3$ are integers $n<y$ where $n$ has no prime factors in $[z,x/3]$ and $n-1$ has no prime factors less than $z$. Let us call the set of such integers $\mathcal{S}$. This is a set which is very similar to the primes, since most elements have a very large prime factor, and it is not clear that there are \textit{any} possible $a_p$ which might cover more than one additional element. Indeed, a typical residue class $a_p\Mod{p}$ for $p>x/3$ will contain \textit{no} elements of $\mathcal{S}$. The problem of showing the existence of unusual residue classes containing many elements of $\mathcal{S}$ leads us to the following toy problem.
\begin{prblm}\label{prblm:ToyProblem}
Given a prime modulus $q$, can we find a residue class $a_q\Mod{q}$ which contains many primes all less than $q(\log{q})^{1/2}$?
\end{prblm}
This toy problem can be answered in the positive by a variant of Theorem \ref{thrm:WeakTuples}. If we choose our linear functions to be of the form $L_i(n)=n+h_iq$ for suitable constants $h_i$, then the existence of a residue class containing many small primes is implied by many of the $L_i$ being simultaneously prime at $n$ for some $n\le q$. The underlying sieve methodology is flexible enough to handle the fact that now the linear functions depend on $q$, and so can solve the toy problem. This correspondingly shows that for each large prime $p$, there are some residue classes $a_p\Mod{p}$ containing many elements of $\mathcal{S}$.\footnote{The work of Green-Tao on linear equations in primes allows one to show for \textit{most} large primes $q$ there are $k$ primes less than $q(\log{q})^{1/2}$ for any fixed $k$, which is sufficient to improve the estimates for large primes.} 

To turn this into an actual covering, we need the residue classes for \textit{different} large primes to be approximately `independent' from one another. To do this we use the probabilistic method, by choosing a residue class at random for each large prime $p\in (x/3,x/2]$, and showing that with high probability this results in approximately independent behavior. Specifically,  we choose $a_p$ randomly with
\[
\mathbb{P}(a_p= a\Mod{p})\propto \sum_{\substack{n\equiv a \Mod{p}\\ L_1(n),\dots,L_k(n)\in \mathcal{S}}}w_{\mathcal{L}_p}(n),
\]
where $w_{L_p}(n)$ is the normalized sieve weight introduced in Section \ref{sec:Sieve} for the functions $L_1,\dots,L_k$ with $L_i(n)=n+h_i p$. We make these choices independently for all $p\in (x/3,x/2]$, so
\begin{align*}
\mathbb{P}\,&(n\text{ not covered by large primes})=\prod_{p\in [x/3,x/2]}\Bigl(1-\sum_{m\equiv n\Mod{p}}w_{\mathcal{L}_p}(m)\Bigr)\\
&\qquad\le \exp\Bigl(-\mathbb{E}\,\#\{p\in [x/3,x/2]:\,n\equiv a_p\Mod{p}\}\Bigr).
\end{align*}
In particular, if the expected number of times any $n\in\mathcal{S}$ is congruent to $a_p\Mod{p}$ for some $p\in (x/3,x/2]$ is at least $t$, then the expected number of elements of $\mathcal{S}$ which are not covered by $a_p\Mod{p}$ for $p\le x/2$ is at most $e^{-t}\#\mathcal{S}$. If $t$ is large, we can then greedily choose residue classes $a_p\Mod{p}$ for $p\in (x/2,x]$ to cover these few remaining elements.

Thus our new strategy for choosing the $a_p$ is:
\begin{enumerate}
\item Choose $a_p=0$ for `medium primes' $p\in [z,x/3]$.
\item Choose $a_p=1$ for `small primes' $p<z$.
\item Choose $a_p$ randomly according to the sieve weights $w_{\mathcal{L}_p}$ independently for `large primes' $p\in(x/3,x/2]$.
\item Choose $a_p$ greedily for `very large primes' $p\in (x/2,x]$.
\end{enumerate}

This gives an arbitrarily large constant improvement over \eqref{eq:RankinBound} provided we show that we can take $t$ arbitrarily large. Calculations similar to those of Section \ref{sec:Sieve} allow us to take $t$ to be of size $ \log{k}$ when we consider $k$ linear functions, and so letting $k$ be large enough we succeed in showing that the constant in \eqref{eq:RankinBound} can indeed be taken to be arbitrarily large.

To get a quantitative improvement over \eqref{eq:RankinBound} we run the essentially the same argument, but we need a version of Theorem \ref{thrm:WeakTuples} which has uniformity with respect to the number $k$ of linear functions we consider as well as uniformity with respect to the coefficients of the linear functions. Such a version of Theorem \ref{thrm:WeakTuples} was established in \cite{MaynardDense}. With some technical modifications, this would ultimately yield a bound
\[
 \sup_{p_n\le X}(p_{n+1}-p_n)\ge c'' \frac{\log{X}\cdot\log\log{X}}{\log\log\log{X}}
\]
for some constant $c''>0$. This is not quite as good as Theorem \ref{thrm:LargeGaps}, because one can improve the quantitative argument further by being more careful about the manner in which we choose the $a_p$ for different large primes. If instead of choosing residue classes independently at random we use ideas based on the `semi-random' or `R\"odl nibble' method from combinatorics, we are able to establish a hypergraph covering lemma which allows our covering by residue classes to have almost no overlaps. After working through the technical details, this ultimately gives an additional improvement of a factor $\log\log\log\log{X}$, and hence gives Theorem \ref{thrm:LargeGaps}.

\section{Limitations}
Both Theorem \ref{thrm:WeakTuples} and Theorem \ref{thrm:LargeGaps} appear to be the qualitative limit of what these methods can give, and require an entirely new approach to do better. Theorem \ref{thrm:SmallGaps} depends on the quantitative aspects of Theorem \ref{thrm:WeakTuples}, and can potentially be improved slightly. New ideas are likely needed to significantly improve upon Theorem \ref{thrm:SmallGaps}, however.

Optimal weights in high-dimensional sieves are poorly understood, and we do not know of general barriers beyond the parity phenomenon. In the context of Theorem \ref{thrm:WeakTuples} the parity phenomenon means that we cannot hope to prove $k/2$ of our linear functions are simultaneously prime based on a sieve argument. In particular, significant new ideas are required to attack the Twin Prime Conjecture. 

Although in principle this leaves open the possibility of having rather more than $\log{k}$ of the linear functions being simultaneously prime, since the Selberg sieve weights seem to perform best in high dimensional sieving situations, we expect that it is unlikely to be possible to much do better than \eqref{eq:NewExpectation} based on weights formed by short divisor sums, even if we cannot prove a direct obstruction. Heuristic arguments show that over fairly general classes of potential Selberg sieve weights, we do not expect to do better than the choice given by \eqref{eq:NewWeights}. Given the choice \eqref{eq:NewWeights}, it is possible to show that a choice of smooth function similar to \eqref{eq:SmoothChoice} is essentially best possible.

For the question of explicit small gaps between primes, there is some potential for further progress. Stronger results about primes in arithmetic progressions allow us to take $R$ in \eqref{eq:NewExpectation} larger, which should reduce the critical value of $k$, and hence the size of the gap. The current record of 246 in Theorem \ref{thrm:SmallGaps} does not make use of the new equidistribution estimates of Zhang or its refinements (but these would likely only lead to small improvements). If we assume optimistic conjectures on the distribution of primes then we can do significantly better. Under the Elliott-Halberstam conjecture one can show gaps of size 12 \cite{Maynard}, and a generalization of the Elliott-Halberstam conjecture to numbers with several prime factors allows us to reach the absolute limit of these methods. Specifically, we have the following.
\begin{thrm}[\cite{Polymath}]\label{thrm:Conditional}
Assume the `Generalized Elliott-Halberstam Conjecture' (see \cite{Polymath}). Then we have
\[
\liminf_n(p_{n+1}-p_n)\le 6.
\]
\end{thrm}
The parity phenomenon makes it impossible for a first moment method of this type to prove a result less than 6, and so Theorem \ref{thrm:Conditional} is the strongest result of this type we can hope to prove along these lines.

All proofs showing the existence of large gaps between primes rely on some variant of Lemma \ref{lmm:CRT}, which allows one to construct a sequence of consecutive composite integers $n_1,\dots,n_1+y$ all with a prime factor of size $O(\log{n_1})$. This places a severe limitation on how large the gaps we produce can be, since we expect that a large gap between primes will involve many composites $n$ whose smallest prime factor is much larger than $\log{n}$. Specifically, Maier and Pomerance \cite{MaierPomerance} conjectured that the largest string of consecutive integers less than $X$ all containing such a small prime factor should be of length
\[
\log{X}\cdot (\log\log{X})^{2+o(1)}.
\]
Therefore we do not expect to be able to produce gaps larger than this without a new approach. The `semi-random' method used in \cite{FGKMT} to show the existence of a good choice residue classes for large primes is essentially as good as one can hope for, so any quantitative progress based on the same overall method would likely require an improvement to Problem \ref{prblm:ToyProblem}, showing the existence of more than $\log\log{q}$ primes less than $q(\log{q})^{1/2}$ in some residue class modulo $q$. A uniform version of Conjecture \ref{cnjctr:KTuples} suggests that there should be residue classes containing roughly $(\log{q})^{1/2}/\log\log{q}$ such primes, but we do not know how to prove this.
\section{Other applications and further reading}

For a more thorough survey of the details of these ideas on small gaps between primes, as well as the ideas behind the breakthrough of Zhang, we refer the reader to the excellent survey articles of Granville \cite{Granville} and Kowalski \cite{Kowalski}. For a more details of the original work of Goldston, Pintz and Y\i ld\i r\i m we refer the reader to survey by Soundararajan \cite{Sound}.

One useful feature of the argument of Section \ref{sec:Sieve} is that the full strength of the Bombieri-Vinogradov theorem was not required to prove bounded gaps between primes. Provided we can estimate primes in arithmetic progressions with modulus of size $x^\epsilon$ on average, we would obtain a version of Theorem \ref{thrm:WeakTuples} with $c\epsilon\log{k}$ of the linear functions are prime. This allows one to show the existence of bounded gaps between primes in many subsets of the primes where one has this type of weaker arithmetic information. A general statement of this type was established in \cite{MaynardDense}. The fact that one can restrict the entire argument to an arithmetic progression also allows one to get some control on the joint distribution of various arithmetic functions. There have been many recent works making use of these flexibilities in the setup of the sieve method, including \cite{Thorner,LemkeOliver,Freiberg,FreibergII,Pollack,Hongze,BakerPollack,MatomakiShao,BakerZhao,ChuaParkSmith,Vatwani,Troupe,PintzII,PintzIII,Huang,BakerZhaoII,BanksFreiberg,BakerFreiberg,Kaptan,Parshall,PollackThompson}. 

New results on long gaps between primes have also found further applications to other situations \cite{BakerFreiberg,MaierRassias,Pratt}. It is hopeful that the ideas behind Theorems \ref{thrm:WeakTuples}-\ref{thrm:ManyGaps} can find further applications in the future.

\section{Acknowledgements}
The author is funded by a Clay Research Fellowship. This article was written whilst the author was a member of the Institute for Advanced Study, Princeton and supported by the National Science Foundation under Grant No. DMS - 1638352.
\bibliographystyle{plain}
\bibliography{ICM}
\end{document}